\documentclass[12pt,onecolumn,twoside]{IEEEtran}

\usepackage{amsmath,amssymb,euscript,psfrag,latexsym,color,graphicx}

\begin{document}
\newtheorem{thm}{Theorem}
\newtheorem{cor}[thm]{Corollary}
\newtheorem{lemma}[thm]{Lemma}
\newtheorem{prop}[thm]{Proposition}
\newtheorem{problem}[thm]{Problem}
\newtheorem{remark}[thm]{Remark}
\newtheorem{defn}[thm]{Definition}
\newtheorem{ex}[thm]{Example}
\newcommand{\ignore}[1]{}
\newcommand{\mR}{{\mathbb R}}
\newcommand{\mN}{{\mathbb N}}
\newcommand{\bR}{{\bf R}}
\newcommand{\rR}{{\rm R}}
\newcommand{\mZ}{{\mathbb Z}}
\newcommand{\mC}{{\mathbb C}}
\newcommand{\fR}{{\mathfrak R}}
\newcommand{\fK}{{\mathfrak K}}
\newcommand{\fG}{{\mathfrak G}}
\newcommand{\fT}{{\mathfrak T}}
\newcommand{\fC}{{\mathfrak C}}
\newcommand{\fF}{{\mathfrak F}}
\newcommand{\cF}{{\mathcal F}}
\newcommand{\fM}{{\mathfrak M}}
\newcommand{\mD}{{\mathbb{D}}}
\newcommand{\cU}{{\mathcal{U}}}
\newcommand{\cL}{{\mathcal{L}}}
\newcommand{\cD}{{\mathcal{D}}}
\newcommand{\cR}{{\mathcal{R}}}
\newcommand{\bJ}{{\mathbb{J}}}
\newcommand{\bI}{{\mathbb{I}}}
\newcommand{\s}{{\rm s}}
\newcommand{\dtheta}{{\frac{d\theta}{2\pi}}}
\newcommand{\ar}{{\rm a}}
\newcommand{\intpi}{{\int_{-\pi}^\pi}}
\newcommand{\cE}{{\mathcal E}}
\newcommand{\me}{{\rm _{ME}}}
\newcommand{\mr}{{\rm _{sm}}}
\def\spacingset#1{\def\baselinestretch{#1}\small\normalsize}

\title{An intrinsic metric\\ for power spectral density functions$^*$}
\author{Tryphon T. Georgiou$^\dag$\thanks{$^\dag$Department of Electrical and Computer Engineering, University of Minnesota, Minneapolis, MN 55455; {\tt tryphon@umn.edu}. $^*$This research has been supported by the NSF and the AFOSR.}}
\date{}
\markboth{August 19, 2006}
{Georgiou: metric on spectral density functions}
\maketitle

\spacingset{1.5}
\begin{abstract} We present an intrinsic metric that quantifies distances between power spectral density functions.
The metric was derived by the author in \cite{Geo2006_metrics} as the geodesic distance between spectral density functions with respect to a particular pseudo-Riemannian metric motivated by a quadratic prediction problem. We provide an independent verification of the metric inequality and discuss certain key properties of the induced topology.
 \end{abstract}

\begin{keywords} Power spectral density functions, intrinsic metric, information geometry.
\end{keywords}



\section{The metric property}\label{sec:introduction}
\PARstart{T}{he} present work builds on a recent report \cite{Geo2006_metrics} where the present author introduced a natural pseudo-Riemannian metric on power spectral density functions of discrete-time stochastic processes, characterized geodesics, and computed geodesic distances. The {\em geodesic distance} between two power spectral density functions $f_i(\theta)$, with $i=1,2$ and $\theta\in[-\pi,\pi]$, was shown to be\\[.02in]
\begin{equation}\label{eq:metric}
d_g(f_1,f_2):=\sqrt{\intpi \left(\log\frac{f_1(\theta)}{f_2(\theta)}\right)^2\dtheta -\left(\intpi \log\frac{f_1(\theta)}{f_2(\theta)}\dtheta\right)^2}.
\end{equation}

\noindent
Below we will provide a direct verification that $d_g(\cdot,\cdot)$ provides a pseudo-metric on the cone of power spectral density functions
\[
\cD:=\{f\,:\, f(\theta)\geq 0 \mbox{ for }\theta\in[-\pi,\pi],\, f\in L_1[-\pi,\pi]\}.
\]
(As usual, $L_1,L_2$ denote Lebeague spaces of integrable and square-integrable functions, respectively.)

The only reason $d_g(\cdot,\cdot)$ is a pseudo-metric and not a metric is because it is insensitive to scaling, i.e., $d_g(f_1,f_2)=d_g(f_1,\kappa f_2)$ for any $\kappa>0$. Thus, it does not differentiate between spectral densities which only differ by a constant nonzero positive factor. Families of spectral density functions related in this way are referred to as spectral rays and form a set
\[
\cR:=\{ f{\rm mod}  \mR_+ \,:\, f(\theta)\geq 0 \mbox{ for }\theta\in[-\pi,\pi],\, f\in L_1[-\pi,\pi]\}
\]
of equivalence classes, and
$d_g(\cdot,\cdot)$ can be used to evaluate distances on $\cR$ via comparing any two representatives on any two given spectral rays. Then, as we will see, $d_g(\cdot,\cdot)$ defines a metric on $\cR$. This metric can be also be readily modified to provide a metric on $\cD$ if for instance, $|\intpi \left(f_1(\theta)-f_2(\theta)\right)\dtheta|$, or the absolute difference of any other generalized means is added on as in
\[d(f_1,f_2):=d_g(f_1,f_2)+|\intpi \left(f_1(\theta)-f_2(\theta)\right)\dtheta|,
\]
to differentiate the effect of scaling.

Before we proceed, we clarify how to evaluate $d_g(\cdot,\cdot)$ on all spectra in $\cD$, including those that may vanish on a subset of the frequency interval $[-\pi,\pi]$ rendering $\log(f_1/f_2)^2$ non-integrable. Clearly, when neither argument of $d_g(f_1,f_2)$ vanishes and $f_i(\theta)$ stays away from zero for $\theta\in[-\pi,\pi]$ and $i=1,2$, then $\log f_i\in L_2[-\pi,\pi]$ and $d_g(f_1,f_2)$ is well defined and finite. But, if either $f_i$ ($i=1,2$) vanishes on $[-\pi,\pi]$ the integrals may diverge. However, since the root-mean-square of {\em any} function, and hence of $\log(f_1/f_2)$ in particular, is always greater than or equal to its arithmetic mean (e.g., see \cite{BB}) it follows that
\begin{equation}\label{eq:LHS}
\sqrt{\intpi \left(\log\frac{f_1(\theta)}{f_2(\theta)}\right)^2\dtheta}\geq \intpi \log\frac{f_1(\theta)}{f_2(\theta)}\dtheta.
\end{equation}
Therefore, (\ref{eq:metric}) gives either a nonnegative real or has to be assigned the value $+\infty$.
In conclusion, we complete the definition of $d_g(\cdot,\cdot)$ as follows: if
\begin{equation}\label{eq:L2condition}
\log \frac{f_1}{f_2}\in L_2[-\pi,\pi],
\end{equation}
in which case the left hand side of (\ref{eq:LHS}) is finite, $d_g(f_1,f_2)$ is evaluated using (\ref{eq:metric}).
If however (\ref{eq:L2condition}) fails then, for consistency with (\ref{eq:metric}), we assign
\begin{equation}\label{eq:infty}
d_g(f_1,f_2):=\infty.
\end{equation}
Clearly, failure of (\ref{eq:L2condition}) can always be traced to at least one of $f_i$ ($i\in\{1,2\}$) failing to satisfy $\log f_i\in L_2[-\pi,\pi]$ (otherwise, necessarily, $\log(f_1/f_2)=(\log f_1 - \log f_2) \in L_2[-\pi,\pi]$).
\vspace*{.1in}

\begin{thm}{\sf $d_g(\cdot,\cdot)$ defines a pseudo-metric on $\cD$ and a metric on $\cR$.
}\end{thm}

\begin{proof} By definition $d(\cdot,\cdot)\in[0,\infty)\cup\{+\infty\}$. It is also easy to observe that
\begin{equation}\label{commutativity}
d_g(f_1,f_2)=d_g(f_2,f_1).
\end{equation}
To see this note that $\log(f_1/f_2)=-\log(f_2/f_1)$ and that (\ref{eq:metric}) is impervious to a sign change in front of the logarithms. Also in case one of $\log(f_1/f_2)$, $\log(f_2/f_1)$ fails to be in $L_2$, so does the other, and again $d_g(f_1,f_2)=d_g(f_2,f_1)$ (both being $\infty$). Thus, (\ref{commutativity}) holds.

When $d_g(f_1,f_2)=0$ the root-mean-square of the function $\log(f_1/f_2)$ is equal to its arithmetic mean, and this only happens (see \cite{BB}) when the function is constant, i.e.,
\begin{eqnarray*}d_g(f_1,f_2)=0&\Rightarrow& \log\frac{f_1(\theta)}{f_2(\theta)} = c \in\mR \mbox{ for all }\theta\in[-\pi,\pi]\\
&\Rightarrow& f_1= \kappa f_2\\
&\Rightarrow& f_1{\rm mod}\mR=f_2{\rm mod}\mR,
\end{eqnarray*}
since $\kappa=e^c$ is a constant. Thus, $d_g(\cdot,\cdot)$ separates the elements of $\cR$.

We finally establish the triangular inequality. So let us consider $f_i\in\cD$ for $i\in\{1,2,3\}$. We will show that
\begin{equation}\label{eq:triangular}
d_g(f_1,f_2)+d_g(f_2,f_3)\geq d_g(f_1,f_3).
\end{equation}
We first argue the case when $d_g(f_1,f_3)=\infty$. It suffices to show that one of the left hand side terms is also infinity.
Assume the contrary, i.e., that
\[\log\frac{f_1}{f_2}\in L_2[-\pi,\pi] \mbox{ as well } \log\frac{f_2}{f_3}\in L_2[-\pi,\pi].
\]
It readily follows that $\log(f_1/f_3)=\log(f_1/f_2)+\log(f_2/f_3)\in L_2[-\pi,\pi]$ which contradicts the assumption that $d_g(f_1,f_3)=\infty$. Thus, at least one of $d_g(f_1,f_2)$, $d_g(f_2,f_3)$ is infinity and the triangular inequality holds. Of course, if $\log(f_1/f_3)$ is finite and any of $d_g(f_1,f_2)$, $d_g(f_2,f_3)$ takes the value $\infty$, the triangular inequality holds anyway.

We now argue the case when all three $d_g(f_1,f_2)$, $d_g(f_2,f_3)$ and $d_g(f_1,f_3)$ are finite.
To this end we square both sides of (\ref{eq:triangular})
and utilize
\begin{equation}\label{eq:logsum}
\log\frac{f_1}{f_3}=\log\frac{f_1}{f_2}+\log\frac{f_2}{f_3},
\end{equation}
to simplify the resulting expression and deduce the following inequality
\begin{eqnarray}\label{eq:triangular2}\nonumber
&&\sqrt{\intpi \left(\log\frac{f_1}{f_2}\right)^2\dtheta - \left(\intpi \log\frac{f_1}{f_2}\dtheta\right)^2}
\sqrt{\intpi \left(\log\frac{f_2}{f_3}\right)^2\dtheta - \left(\intpi \log\frac{f_2}{f_3}\dtheta\right)^2}\\
&\geq&\intpi \left(\log\frac{f_1}{f_2}\log\frac{f_2}{f_3}\right)\dtheta - \left(\int \log\frac{f_1}{f_2}\dtheta\right)\left(\int \log\frac{f_2}{f_3}\dtheta\right).
\label{questionmark}
\end{eqnarray}
Thus the two inequalities (\ref{eq:triangular2}) and (\ref{eq:triangular}) are equivalent to one another, and therefore, in order to ascertain (\ref{eq:triangular}) it suffices to establish the validity of (\ref{eq:triangular2}).

To this end, let $\alpha:=\log(f_1/f_2)$, $\beta:=\log(f_2/f_3)$ and rewrite (\ref{eq:triangular2}) in the form
\begin{eqnarray}\nonumber
&&\sqrt{\intpi \alpha^2\dtheta-(\intpi \alpha \dtheta)^2}
\sqrt{\intpi \beta^2\dtheta-(\intpi \beta \dtheta)^2}\\
&\geq &
\intpi \left(\alpha\beta\right) \dtheta -\left(\intpi \alpha \dtheta\right)\left( \intpi \beta \dtheta\right).
\label{eq:generalinequality}
\end{eqnarray}
Since (\ref{eq:generalinequality}) is homogeneous in both $\alpha$ and $\beta$, scaling of either leaves it unaffected. Therefore, if
\begin{eqnarray*}\sigma_\alpha&=&\sqrt{\intpi \alpha^2\dtheta-(\intpi \alpha \dtheta)^2}\\
\sigma_\beta&=&\sqrt{\intpi \beta^2\dtheta-(\intpi \beta \dtheta)^2}
\end{eqnarray*}
and $a:=\frac{1}{\sigma_\alpha}\alpha$, $b:=\frac{1}{\sigma_\beta}\beta$,
the inequality (\ref{eq:generalinequality}) is equivalent to
\begin{eqnarray}1
&\geq &
\intpi \left(ab\right) \dtheta -\left(\intpi a \dtheta\right)\left( \intpi b \dtheta\right)
\label{eq:generalinequality2}
\end{eqnarray}
with the side conditions
\begin{eqnarray}\label{sigmaa}
\intpi a^2\dtheta-(\intpi a \dtheta)^2&=&1, \mbox{ and}\\
\intpi b^2\dtheta-(\intpi b \dtheta)^2&=&1\label{sigmab}.
\end{eqnarray}
But the validity of (\ref{eq:generalinequality2}) follows trivially from the standard inequality
\[\intpi (a-b)^2\dtheta\geq \left(\intpi (a-b)\dtheta\right)^2\]
after we expand the squares on both sides and use (\ref{sigmaa}) and (\ref{sigmab}) to simplify the resulting expressions.
Thus (\ref{eq:generalinequality2}) with (\ref{sigmaa}-\ref{sigmab}) holds $\Rightarrow$ (\ref{eq:generalinequality}) $\Rightarrow$ (\ref{questionmark}) $\Rightarrow$ (\ref{eq:triangular}), and this completes the proof.
\end{proof}

{\em Remark:} The definition of $d_g(\cdot,\cdot)$ distinguishes two classes of power spectral densities according to whether their logarithm is square integrable or not.
The first class,
$\cD_{\rm interior}:=\{f\in\cD\;:\;\log f\in L_2[-\pi,\pi]\}$,
can be thought of as ``interior'' points lying to within a finite distance from one another, and to within a finite distance from constant non-zero power spectral densities. The second class, with logarithms that fail to be square integrable, contains power spectral densities which lie at an infinite distance from any density in $\cD_{\rm interior}$.
On the other hand, power spectal densities are traditionally differentiated according to whether the underlying process is deterministic or not. More specifically, a stochastic process is said to be {\em non-deterministic} (in the sense of Kolmogoroff) if the
variance of the one-step-ahead prediction error cannot be made arbitrarily small.
In turn, this property is characterized by the log-integrability of the corresponding power spectral density function (see \cite{GrenanderSzego,StoicaMoses}), i.e.,
$\cD_{\rm non-deterministic}:=\{f\in\cD\;:\; \log f\in L_1[-\pi,\pi]\}$.
Thus, it is interesting to observe that
$\cD_{{\rm interior}}\subset \cD_{\rm non-deterministic}$ and hence, finite neighborhoods of elements in $\cD_{{\rm interior}}$ contain non-deterministic power spectra {\em only}.

\section{Riemannian geometry, geodesics, and intrinsic metrics}
We now explain the geometric significance of $d_g(\cdot,\cdot)$
recapitulating some of the development in \cite{Geo2006_metrics}. The starting point that led to (\ref{eq:metric}) is a prediction problem and the degradation of the variance of the prediction error when the design of the predictor is based on the wrong choice among two alternatives. More specifically, let $f_1,f_2$ represent spectral density functions of discrete-time zero-mean stochastic processes $u_{f_i}(k)$ ($i\in\{1,2\}$ and $k\in\mZ$), and let
$p_{f_i}(\ell)$ ($\ell\in\{1,2,3,\ldots\}$) represent values for the coefficients that minimize the linear prediction error variance
\[\cE\{|u_{f_i}(0)-\sum_{\ell=1}^\infty p(\ell) u_{f_i}(-\ell)|^2\}.
\]
Thus, the optimal set of coefficients depends on the power spectral density function of the process, a fact which is duly acknowledged by the subscript in the notation $p_{f_i}(\ell)$.
Here, as usual, $\cE\{\;\}$ denotes the expectation operator.
It is reasonable to consider as a distance between $f_1$ and $f_2$ the degradation of predictive error variance when the coefficients $p(\ell)$ are selected assuming one of the two, and then used to predict a stochastic process corresponding to the other spectral density function. The ratio of the ``degraded'' predictive error variance over the optimal error variance
\[
\rho(f_1,f_2):=\frac{\cE\{u_{f_1}(0)-\sum_{\ell=1}^\infty p_{f_2}(\ell) u_{f_1}(-\ell)|^2\}}
       {\cE\{u_{f_1}(0)-\sum_{\ell=1}^\infty p_{f_1}(\ell) u_{f_1}(-\ell)|^2\}}
\]
turns out to be equal to the ratio of the arithmetic over the geometric means of the fraction of the two spectral density functions, namely
\[
\rho(f_1,f_2)=\frac{\intpi \frac{f_1(\theta)}{f_2(\theta)}\dtheta}
                                            {\exp\left(\intpi \log (\frac{f_1(\theta)}{f_2(\theta)})\dtheta\right)},
\]
see \cite{Geo2006_metrics}.

The logarithm $\log\rho(f_1,f_2)=:\delta_{\rm a/g}(f_1,f_2)$ (where the subscript ${\rm a/g}$ signifies arithmetic/geometric) represents a measure of dissimilarity between the ``shapes'' of $f_1$ and $f_2$ and, can be viewed, as analogous to ``divergences'' of Information Theory. Indeed,
\[
\delta_{\rm a/g}(f_1,f_2)=\log\left(\intpi \frac{f_1(\theta)}{f_2(\theta)}\dtheta\right)-
                                            \intpi \log \left(\frac{f_1(\theta)}{f_2(\theta)}\right)\dtheta
\]
vanishes only when $f_1/f_2$ is constant on $[-\pi,\pi]$ and is positive otherwise.
Considering the distance $\delta_{\rm a/g}(f,f+\Delta)$ between a nominal power spectral density $f$ and a perturbations $f+\Delta$, and eliminating cubic terms and beyond, leads (modulo a scaling factor of $2$) to the Riemannian pseudo-metric in $\cD_{\rm interior}$ which is given by the following quadratic differential form
\begin{equation}\label{Gmetric}
g_f(\Delta):= \intpi \left(\frac{\Delta(\theta)}{f(\theta)}\right)^2\dtheta -\left(\intpi \frac{\Delta(\theta)}{f(\theta)}\dtheta\right)^2.
\end{equation}
Interestingly, geodesic paths $f_\tau$ ($\tau\in[0,1]$) connecting spectral densities $f_0,f_1$ and having minimal length
\[\sqrt{2}\int_0^1\sqrt{\delta_{\rm a/g}(f_\tau,f_{\tau+d\tau})}=\int_0^1 \sqrt{g_{f_\tau}(\frac{\partial f_\tau}{\partial \tau})}d\tau
\]
can be explicitely computed \cite{Geo2006_metrics}. They turn out to be logarithmic intervals
\begin{equation}\label{eq:loginterval}
f_\tau(\theta)=f_0^{1-\tau}(\theta)f_1^\tau(\theta) \mbox{ for }\tau\in[0,1],
\end{equation}
between the two extreme points. Furthermore, the length along such geodesics is precisely $d_g(f_0,f_1)$ as given in (\ref{eq:metric}).

The closed form of the geodesic path allows us to verify directly that
any two power spectral densities $f_0,f_1$, at a finite distance from one another, can be connected with a path of the same length.
A topological space with such a property is said to be a {\em length-space} and the
metric is said to be {\em intrinsic}.
The fact that $d_g(\cdot,\cdot)$ is intrinsic can be readily verified and this is done below.

\begin{prop}{\sf $d_g(\cdot,\cdot)$ is intrinsic on $\cD$ and $\cR$.
}\end{prop}

\begin{proof}
By direct substitution into (\ref{eq:metric})
we can verify that for any $f_0,f_1$ such that $d_g(f_0,f_1)<\infty$, any $\tau\in[0,1]$, and with $f_\tau$ defined as in (\ref{eq:loginterval}),
$d_g(f_0,f_\tau)=\tau d_g(f_0,f_1)$, $d_g(f_\tau,f_1)=(1-\tau) d_g(f_0,f_1)$, and even
$d_g(f_\tau,f_{\tau+d\tau})=d_g(f_0,f_1)d\tau$.
It readily follows that the length of the path $\int_0^1 d_g(f_\tau,f_{\tau+d\tau})=\int_0^1 d_g(f_0,f_1)d\tau$ equals the distance $d_g(f_0,f_1)$ between the end points.
\end{proof}

{\em Remark:}
Besides $\delta_{\rm a/g}(f_1,f_2)$, several other ``divergences'' have been introduced
in \cite{Geo2006_metrics} as appropriate distance measures (though not metrics). First the symmetrized version
\[\delta(f_1,f_2)=\delta_{\rm a/g}(f_1,f_2)+\delta_{\rm a/g}(f_2,f_1),\]
and then,
\[
\delta_{r,s}(f_1,f_2):= \log\sqrt[r] {\intpi \left(\frac{f_1}{f_2}\right)^r\dtheta} -
\log \sqrt[s]{\intpi \left(\frac{f_1}{f_2}\right)^s\dtheta}
\]
involving comparison of other generalized means.
It is interesting to point out that
the quadratic terms of  $\delta_{r,s}(f,f+\Delta)$, $\delta_{\rm a/g}(f,f+\Delta)$ and $\delta(f,f+\Delta)$ in the ``perturbation'' $\Delta$, all turn out to be identical (modulo a scaling). Hence, they all lead to the same Riemannian pseudo-metric (\ref{Gmetric}).

\section{Concluding thoughts}
It is interesting to compare the differential structure on power spectral density functions which we have introduced above with the corresponding differential structure of ``Information Geometry.'' In Information Geometry $f(\theta)$ corresponds to a probability density on $[-\pi,\pi]$
and the natural Riemannian metric is the {\em Fisher information metric} is (cf.\ \cite[page 28]{amari}) which is this case would be
\begin{eqnarray}\nonumber
g_{{\rm Fisher},f}(\Delta)&=&\intpi \left(\frac{\Delta(\theta)}{f(\theta)}\right)^2f(\theta)\dtheta\\
&=&\intpi\frac{\Delta(\theta)^2}{f(\theta)}\dtheta\label{Fmetric}
\end{eqnarray}
(with $\intpi f(\theta)\dtheta =1$ and $\intpi \Delta(\theta)\dtheta=0$ since both $f$, $f+\Delta$ need to be probability densities). Direct comparison reveals that the powers of $f(\theta)$ in (\ref{Gmetric}) and (\ref{Fmetric}) are different. Thus, it is curious and worth underscoring that in either differential structure, geodesics and geodesic lengths can be computed.

\end{document}